%% file: draft.tex
\documentclass[12pt]{article}
\usepackage{epsfig}
\usepackage{latexsym}
\newcommand{\defeq}{\stackrel{\rm def}{=}}
\newcommand{\dstyle}[1]{\displaystyle{#1}}
\newcommand{\e}{{\rm e}}
\newcommand{\del}[1]{\Delta #1}
\newcommand{\bsquare}{\hbox{\rule{6pt}{6pt}}}
\newcommand{\lto}{\leftarrow}
\begin{document}
\begin{center}
EASY ESTIMATION BY A NEW PARAMETERIZATION
FOR THE THREE-PARAMETER LOGNORMAL DISTRIBUTION
\end{center}
\begin{center}
YOSHIO KOMORI and HIDEO HIROSE\\
Department of Control Engineering and Science,\\
Kyushu Institute of Technology,
Iizuka 820-8502, Japan
\end{center}
%

A new parameterization and algorithm are proposed for seeking the primary
relative maximum of the likelihood function in the three-parameter
lognormal distribution. The parameterization yields the dimension reduction
of the three-parameter estimation problem to a two-parameter estimation
problem on the basis of an extended lognormal distribution. The algorithm
provides the way of seeking the profile of an object function in the
two-parameter estimation problem. It is simple and numerically stable
because it is constructed on the basis of the bisection method. The profile
clearly and easily shows whether a primary relative maximum exists or not,
and also gives a primary relative maximum certainly if it exists.
%
\vskip 4mm

\vskip 3mm

\noindent {\bf 1   INTRODUCTION}

The three-parameter lognormal distribution is one of the most
important distributions
in many fields. With 
a variable $x$ and three parameters $\alpha$, $\beta$ and $\gamma$, the
probability density function is expressed by
$$
f(x;\alpha, \beta, \gamma) \defeq \frac{1}{\sqrt{2\pi}(x-\alpha)\beta}
\exp\left[-\frac{\left\{\ln{\bigl((x-\alpha)/\gamma\bigr)}\right\}^{2}}
{2\beta^{2}}\right],\quad x>\alpha,\ \beta>0,\ \gamma>0
\eqno(1.1)
$$
and the likelihood function is expressed by $L(\alpha,\beta,\gamma)
\defeq\prod_{i=1}^{n}f(x_{i};\alpha, \beta, \gamma)
$. Here, $x_{i}\ (1\leq i\leq n)$ stand for independent
observations. Without loss of 
generality, we assume $x_{1}>x_{2}\geq\cdots\geq x_{n-1}>x_{n}$.

Since $\ln(X-\alpha)$ obeys a normal distribution if a random
variable $X$ obeys a lognormal distribution, $L(\alpha,\beta,\gamma)$
achieves its maximum at a point $\bigl(\alpha_{0},\hat{\beta}(\alpha_{0}),
\hat{\gamma}(\alpha_{0})\bigr)$ provided that $\alpha$ is fixed
to $\alpha_{0}$, where
$$
\hat{\beta}(\alpha)\defeq\sqrt{\frac{1}{n}
\sum_{i=1}^{n}\bigl\{\ln(x_{i}-\alpha)-\ln\hat{\gamma}
(\alpha)\bigr\}^{2}}\quad \mbox{and}\quad 
\hat{\gamma}(\alpha)\defeq \exp\left[\frac{1}{n}
\sum_{i=1}^{n}\ln(x_{i}-\alpha)\right].
$$
Consequently, if we want to obtain the maximum likelihood estimate,
it suffices to find an $\alpha$ such that $\hat{L}(\alpha)\defeq L\bigl(
\alpha,\hat{\beta}(\alpha),\hat{\gamma}(\alpha)\bigl)$ achieves its
maximum. However, because that $\hat{L}(\alpha)$ 
 $\to$ $\infty$ as $\alpha$ $\to$ $x_{n}-0$, $L(\alpha,\beta,\gamma)$
becomes unbounded. Furthermore, the other parameters then lead to
inadmissible values.

Hill (1963), using the Bayes theorem, has
given a statistical implication
of $\bigl(\hat{\alpha},\hat{\beta}(\hat{\alpha}),
\hat{\gamma}(\hat{\alpha})
\bigr)$ at which $L(\alpha,\beta,\gamma)$ has its maximum in the region 
except the singular region; for a small
$\delta>0$, at $\alpha=\hat{\alpha}$ has $\hat{L}(\alpha)$ its relative and
absolute maximum under the condition $x_{n}-\alpha>\delta$. The
point $\bigl(\hat{\alpha},\hat{\beta}(\hat{\alpha})
,\hat{\gamma}(\hat{\alpha})\bigr)$ is used instead of the
maximum likelihood estimate, and it is called the primary relative
maximum (PRM) or the local maximum likelihood estimate
of the likelihood function.

Displaying $\hat{L}(\alpha)$ is an 
effective way for finding $\hat{\alpha}$, but the search may be difficult 
because that the shape of $\hat{L}(\alpha)$ is
complicated in some cases
depending on data
sets (Cheng and Iles, 1990; Hill, 1963; Johnson, Kotz and Balakrishnan, 1994).
For example, $\hat{L}(\alpha)$ for Data 3 (Table I) increases
very slowly as $\alpha$ becomes small, furthermore it begins
to oscillate when $\alpha$ becomes sufficiently small due to numerical
errors (Fig. 1). On the other hand, $\hat{L}(\alpha)$ for Data 4
(Table 1)
attains its maximum at a point which is closely
near to the end-point of the domain of definition
(Fig. 2). In this case it is necessary
to magnify carefully the graph to a big scale not to miss its
maximum. Besides, if an iterative solver like Newton's method
is used to find the maximum, some difficulties can happen.
This is also one of the reasons why many researchers
tackled this estimation problem.

%

For problems to seek $\hat{\alpha}$, that is, the one-parameter estimation
problems for the lognormal
distribution, Wingo (1975, 1976, 1984) has proposed
a computing method to avoid the singular
range $x_{n}-\alpha\leq\delta$ by adopting a penalty
function. On the other hand, for 
problems to seek simultaneously
$\hat{\alpha}$, $\hat{\beta}$
and $\hat{\gamma}$, where $(\hat{\alpha},\hat{\beta},
\hat{\gamma})$ is a
PRM, that is, the three-parameter estimation
problems (Lambert, 1964) for the distribution, Munro and
Wixley (1970) have proposed a 
parameterization to improve the convergency of many iterative
methods (Eastham, LaRiccia and Schuenemeyer, 1987; Hirose, 1997).
As seen now, there are two ways 
for dealing with the parameter estimation of the distribution for
complete data. Besides 
these, Giesbrecht and Kempthorne (1976) have proposed replacing complete
data with grouped data  to avoid the singularity described above.

The use of Munro and Wixley's parameterization, that is, the substitutions of
$\alpha=\mu-\sigma/\lambda$, $\beta = \lambda$
and $\gamma=\sigma/\lambda$ into (1.1) yield
$$
f(x;\mu-\sigma/\lambda,\lambda,\sigma/\lambda)
=
\frac{1}{\sqrt{2\pi}\bigl\{\sigma+\lambda(x-\mu)\bigr\}}
\exp\left[
-\frac{\left\{\ln\bigl(\sigma+\lambda(x-\mu)\bigr)-\ln\sigma\right\}^{2}}
{2\lambda^{2}}
\right].
\eqno(1.2)
$$
This can be extended by allowing $x<\mu-\sigma/\lambda$, then the
generalization permits $\lambda$ to be negative. In this way, we obtain
the density function for the extended lognormal 
distribution permitting that $\lambda\neq 0$ and $\sigma>0$. Let
$\tilde{f}(x;\lambda,\mu,\sigma)$ be this density function
and $\tilde{L}(\lambda,\mu,\sigma)$ the likelihood function. Cheng
and Iles (1990) have shown that
as $\lambda\to 0$, $\tilde{f}(x;\lambda,\mu,\sigma)$ leads
to the normal distribution with mean $\mu$ and variance
$\sigma^{2}$, which is called the embedded distribution. They
have also investigated tests of statistical hypothesis
to see whether the embedded model should be used.
Consequently, Munro and Wixley's parameterization can not only improve the
convergency of many iterative methods but also cope with the
embedded problem.

In the present article we propose a new reparameterization
of the extended lognormal distribution to change
the three-parameter estimation problem to a two-parameter
estimation problem. Because that the reparameterization also permits $\lambda$
to be estimated negatively, it can cope with even data that cause the
embedded problem. In addition, on the two-parameter estimation
problem we propose an algorithm to obtain stably the profile of an object
function. This makes it possible to seek certainly a PRM if it
exists or to show clearly it does not exist.

\begin{figure}[t]
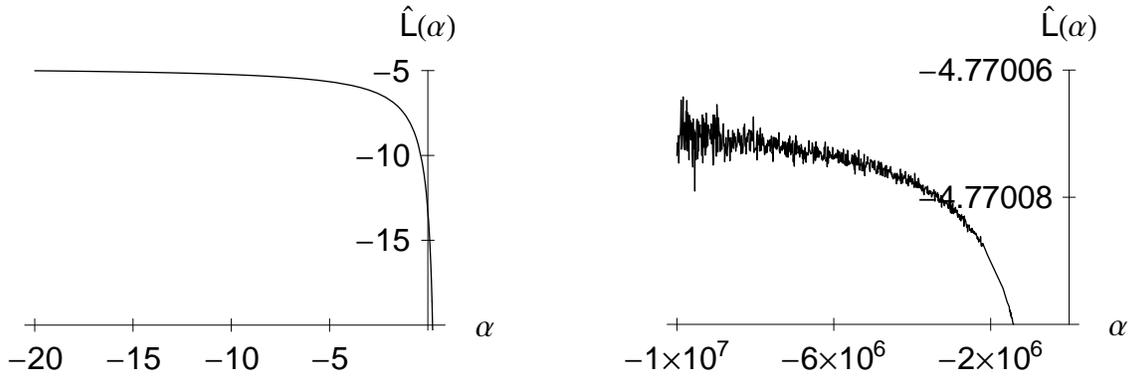

\unitlength 5mm
\begin{center}
\begin{picture}(32,10)(0,0)
\put(0.8,0){\epsfig{file=Smith_figL.epsi,scale=1}}
\put(17.2,0){\epsfig{file=Smith_figS.epsi,scale=1}}
\end{picture}
\caption{$\hat{L}(\alpha)$ for Data 3.}
\end{center}
\end{figure}

\begin{figure}[t]
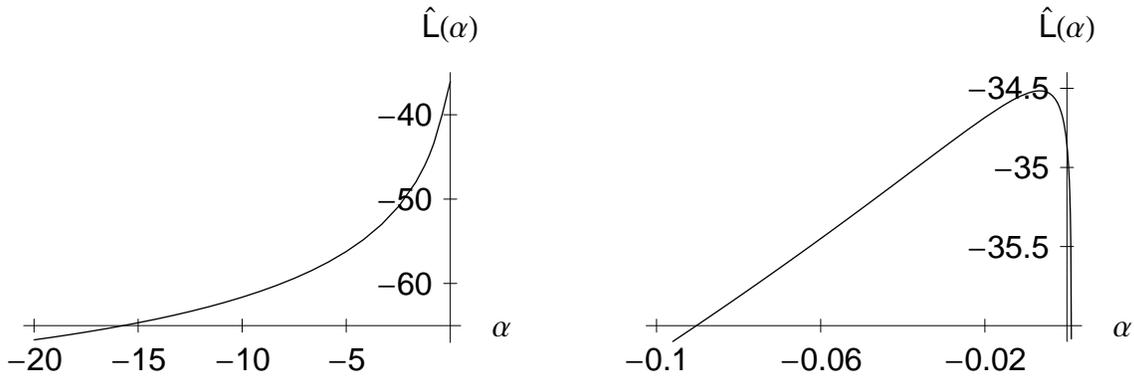

\unitlength 5mm
\begin{center}
\begin{picture}(32,10)(0,0)
\put(0.8,0){\epsfig{file=Menon_figL.epsi,scale=1}}
\put(17.2,0){\epsfig{file=Menon_figS.epsi,scale=1}}
\end{picture}
\caption{$\hat{L}(\alpha)$ for Data 4.}
\end{center}
\end{figure}

In Section 2 we describe the reparameterization and an object function
to be maximized, and give a theorem, which is useful to construct the
algorithm. In Section 3 we introduce the algorithm to obtain the profile
of the object function and a PRM. In Section 4
we challenge some estimation examples and perform Monte Carlo simulation
experiments. A discussion and summary are given in the last
two sections.
\vskip 3mm

\noindent {\bf 2 TWO-PARAMETER ESTIMATION}

In this section we analyze a function maximized to find a PRM. First
of all we introduce the function.

Set $\tau =\sigma-\lambda\mu$ and $s = \ln\sigma$, and
define $\bar{f}(x;\lambda,\tau,s)\defeq
\tilde{f}(x;\lambda,(\e^{s}-\tau)/\lambda,\e^{s})$:
$$
\bar{f}(x;\lambda,\tau,s)\defeq
\frac{1}{\sqrt{2\pi}(\lambda x+\tau)}
\exp\left[
-\frac{\bigl\{\ln(\lambda x+\tau)-s\bigr\}^{2}}
{2\lambda^{2}}
\right],\qquad \lambda\neq 0.
$$
By arranging $\ln\bar{L}(\lambda,\tau,s)
\defeq\sum_{i=1}^{n}\ln\bar{f}(x_{i};\lambda,\tau,s)$, we obtain
\begin{eqnarray*}
\ln\bar{L}(\lambda,\tau,s)&=&
-\frac{n}{2\lambda^{2}}
\left\{s-\frac{1}{n}\sum_{i=1}^{n}\ln(\lambda x_{i}+\tau)\right\}^{2}
-n\ln\sqrt{2\pi}\nonumber\\
&&\makebox[1em]{}+\frac{1}{2n\lambda^{2}}
\left\{\sum_{i=1}^{n}\ln(\lambda x_{i}+\tau)\right\}^{2}
-\frac{1}{2\lambda^{2}}
\sum_{i=1}^{n}\bigl\{\ln(\lambda x_{i}+\tau)\bigr\}^{2}
-\sum_{i=1}^{n}\ln(\lambda x_{i}+\tau).
\end{eqnarray*}
Only the first term depends on $s$ in the right-hand
side of the above equation. And this term has the maximum value $0$
when $s=(1/n)\sum_{i=1}^{n}\ln(\lambda x_{i}+\tau)$. Hence it
suffices to maximize the sum of the third, the fourth and the fifth
terms in the equation. Expressing the sum
by $F(\lambda,\tau)$, let us deal with it:
$$
F(\lambda,\tau)\defeq
\frac{1}{2n\lambda^{2}}
\left\{\sum_{i=1}^{n}\ln(\lambda x_{i}+\tau)\right\}^{2}
-\frac{1}{2\lambda^{2}}
\sum_{i=1}^{n}\bigl\{\ln(\lambda x_{i}+\tau)\bigr\}^{2}
-\sum_{i=1}^{n}\ln(\lambda x_{i}+\tau).
$$

Next, we introduce a useful theorem for constructing
an algorithm searching for a PRM

\noindent {\bf Theorem} \ We set $\bar{x}\defeq(1/n)\sum_{i=1}^{n}x_{i}$, and
define $\tau_{U}^{+}(\lambda)$ and $\tau_{U}^{-}(\lambda)$ as follows:
$$
\tau_{U}^{+}(\lambda)\defeq -\lambda x_{n}
\left(\frac{1-\bar{x}/x_{n}\e^{-\lambda^2}}
{1-\e^{-\lambda^2}}\right)\ \mbox{for}\ \lambda>0,
\quad\tau_{U}^{-}(\lambda)\defeq -\lambda x_{1}
\left(\frac{1-\bar{x}/x_{1}\e^{-\lambda^2}}
{1-\e^{-\lambda^2}}\right)\ \mbox{for}\ \lambda<0.
$$
Then, the following statements hold.
\begin{enumerate}
\renewcommand{\labelenumi}{\arabic{enumi}$)$}
\item
$
\dstyle{
\frac{\partial F}{\partial \tau}
\bigl(\lambda,\tau_{U}^{+}(\lambda)\bigr)<0
\quad \mbox{for $\lambda>0$},\qquad
\frac{\partial F}{\partial \tau}
\bigl(\lambda,\tau_{U}^{-}(\lambda)\bigr)<0
\quad \mbox{for $\lambda<0$}.
}
$
\item
$
\dstyle{
\lim_{\tau\to-\lambda x_{n}+0}
\frac{\partial F}{\partial \tau}(\lambda,\tau)=+\infty
\ \mbox{for $\lambda>0$},\quad
\lim_{\tau\to-\lambda x_{1}+0}
\frac{\partial F}{\partial \tau}(\lambda,\tau)=+\infty
\ \mbox{for $\lambda<0$}.
}
$
\item
$
\dstyle{
\frac{\partial F}{\partial \tau}(\lambda,\tau)<0
}
$ for any point $(\lambda,\tau)$ that satisfies
$\dstyle{\frac{\partial^{2}F}{\partial \tau^{2}}(\lambda,\tau)=0}$.
\item
$
\dstyle{
\lim_{\lambda\to+\infty}F\bigl(\lambda,\tau_{U}^{+}(\lambda)\bigr)=+\infty,
\qquad
\lim_{\lambda\to-\infty}F\bigl(\lambda,\tau_{U}^{-}(\lambda)\bigr)=+\infty.
}
$
\item
$\dstyle{F_{0}(\tau)\defeq\lim_{\lambda\to\pm 0}F(\lambda,\tau)}$
achieves the relative maximum when
$$
\tau = \tau^{\ast}\defeq\frac{1}{n}
\sqrt{\sum_{i=1}^{n-1}\sum_{j=i+1}^{n}(x_{i}-x_{j})^{2}}.
$$
\end{enumerate}
\hfill \bsquare

\noindent {\it Proof.} For $\lambda>0$, we can obtain
$$
\frac{\partial F}{\partial \tau}(\lambda,\tau)<
\frac{1}{\lambda^{2}}\left(
\sum_{i=1}^{n}\frac{1}{\lambda x_{i}+\tau}
\right)
\left(
\ln\frac{\lambda\bar{x}+\tau}{\lambda x_{n}+\tau}-\lambda^{2}
\right)
\eqno(2.1)
$$
by ($A$.1), Jensen's inequality and
$$
\sum_{i=1}^{n}\frac{1}{\lambda x_{i}+\tau}
\left(
\ln \frac{\lambda\bar{x}+\tau}{\lambda x_{i}+\tau}-\lambda^{2}
\right)<
\left(
\sum_{i=1}^{n}\frac{1}{\lambda x_{i}+\tau}
\right)
\left(
\ln\frac{\lambda\bar{x}+\tau}{\lambda x_{n}+\tau}-\lambda^{2}
\right).
$$
By setting the  right-hand side of (2.1) equals 0 and arranging
it, we can see $\tau=\tau_{U}^{+}(\lambda)$. This
leads to the inequality in 1). The proof for $\lambda<0$ is similar.

The statements in 2), 4) and 5) are obtained by direct calculations.

When $\dstyle{\frac{\partial^2 F}{\partial \tau^2}(\lambda,\tau)=0}$, ($A$.2)
in Appendix is equivalent to
$$
\frac{1}{n\lambda^2}
\left\{
\sum_{i=1}^{n}\ln(\lambda x_{i}+\tau)
\right\}=
\frac{1}{\lambda^2}
\left[
\frac{\dstyle{\left(\sum_{i=1}^{n}\frac{1}{\lambda x_{i}+\tau}\right)^2}}
     {\dstyle{n\sum_{i=1}^{n}\frac{1}{(\lambda x_{i}+\tau)^2}}}
+
\frac{\dstyle{\sum_{i=1}^{n}\frac{1}{(\lambda x_{i}+\tau)^2}
      \ln(\lambda x_{i}+\tau)}}
     {\dstyle{\sum_{i=1}^{n}\frac{1}{(\lambda x_{i}+\tau)^2}}}
-1
\right]+1.
$$
The substitution of this into ($A$.1) yields
\begin{eqnarray*}
\lefteqn{
\frac{\partial F}{\partial \tau}(\lambda,\tau)=
\frac{1}{\lambda^2}
\left(\sum_{i=1}^{n}\frac{1}{\lambda x_{i}+\tau}\right)
\left[
\frac{n}{\dstyle{\sum_{i=1}^{n}\frac{1}{(\lambda x_{i}+\tau)^2}}}
\left\{
\left(\frac{1}{n}\sum_{i=1}^{n}\frac{1}{\lambda x_{i}+\tau}\right)^2
-\frac{1}{n}\sum_{i=1}^{n}\frac{1}{(\lambda x_{i}+\tau)^2}
\right\}
\right.
}\\
&&\makebox[12em]{}
\left.
-\frac{\dstyle{\sum_{i=1}^{n}\frac{1}{(\lambda x_{i}+\tau)^2}
                            \ln\frac{1}{\lambda x_{i}+\tau}}}
      {\dstyle{\sum_{i=1}^{n}\frac{1}{(\lambda x_{i}+\tau)^2}}}
+\frac{\dstyle{\sum_{i=1}^{n}\frac{1}{(\lambda x_{i}+\tau)}
                            \ln\frac{1}{\lambda x_{i}+\tau}}}
      {\dstyle{\sum_{i=1}^{n}\frac{1}{(\lambda x_{i}+\tau)}}}
\right].
\end{eqnarray*}
In the bracket, the first term is negative. Furthermore, we can
see the sum of the second and third ones is also negative by
noting that $\dstyle{\frac{\sum_{i} a_{i}\ln a_{i}}{\sum_{i} a_{i}}}
\leq\dstyle{\frac{\sum_{i} a_{i}^{2}\ln a_{i}}{\sum_{i} a_{i}^{2}}}$ holds
for any $a_{i}(>0)$, with equality if and only if all $a_{i}$'s are equal.
\hfill $\Box$

The statements 1) and 2) indicate that for each $\lambda>0$ there exists
a solution, say $\tau_{0}(\lambda)$, of $\frac{\partial F}{\partial \tau}
(\lambda,\tau)=0$ in $\bigl(-\lambda x_{n},\tau_{U}^{+}(\lambda)\bigr)$,
and for each $\lambda<0$ it exists in $\bigl(-\lambda x_{1},
\tau_{U}^{-}(\lambda)\bigr)$.
On the other hand, $\frac{\partial F}{\partial \tau}
(\lambda,\tau)<0$  for each $\lambda>0$
and any $\tau\geq\tau_{U}^{+}(\lambda)$
because of 1) and (2.1), and
it holds for each $\lambda<0$
and any $\tau\geq\tau_{U}^{-}(\lambda)$ 
because of similar reasons.
Furthermore, $\frac{\partial F}{\partial \tau}(\lambda,\tau)<0$
holds for each $\lambda>0$
and $\tau\in(\tau_{0}(\lambda),\tau_{U}^{+}(\lambda))$
and for each $\lambda<0$
and $\tau\in(\tau_{0}(\lambda),\tau_{U}^{-}(\lambda))$
since 1) and 3). Thus, $\tau_{0}(\lambda)$ is the unique
solution of $\frac{\partial F}{\partial \tau}(\lambda,\tau)=0$ for
each $\lambda$. In addition, from these facts, 1) and 4),
$\dstyle{\lim_{\lambda\to\pm\infty}F(\lambda,\tau_{0}(\lambda))=
\infty}$. 
The statement 5) will
be used at the beginning in the algorithm stated below.
Finally, note that
the intervals in which
a $\tau_{0}(\lambda)$ exits, that is,
$(-\lambda x_{n},\tau_{U}^{+}(\lambda))$ and
$(-\lambda x_{1},\tau_{U}^{-}(\lambda))$ become
rapidly narrower as $|\lambda|$ becomes larger.
For instance, when $\lambda=6$, the width of
$(-\lambda x_{n},\tau_{U}^{+}(\lambda))$ is
$\frac{6}{\e^{36}-1}(\bar{x}-x_{n})$.

\clearpage

\noindent {\bf 3  AN ALGORITHM FOR SEEKING THE PROFILE OF $F$}

The theorem can be used to seek the profile
of $F(\lambda,\tau)$ concerning $\lambda$ with the bisection method.
In $\lambda>0$, the procedure is written as follows:
\begin{enumerate}
\renewcommand{\labelenumi}{\arabic{enumi}$)$}
\item
$\tau\lto\delta_{0}$, $\lambda\lto\varepsilon_{0}>0$.
\item
If $\lambda>\lambda_{max}^{+}$, end. Otherwise,
$\tau_{min}\lto -\lambda x_{n}$, $\tau_{max}\lto \tau_{U}^{+}(\lambda)$.
\item
If $\tau_{min}<\tau<\tau_{max}$ and
$\dstyle{\frac{\partial F }{\partial\tau}(\lambda,\tau)> 0}$,
$\tau_{min}\lto\tau$.
If $\tau_{min}<\tau<\tau_{max}$ and
$\dstyle{\frac{\partial F }{\partial\tau}(\lambda,\tau)\leq 0}$,
$\tau_{max}\lto\tau$.
\item
If $\dstyle{\frac{\partial F }{\partial\tau}
\bigl(\lambda,(\tau_{min}+\tau_{max})/2\bigr)}>0$,
then $\tau_{min}\lto(\tau_{min}+\tau_{max})/2$.
Otherwise, $\tau_{max}\lto(\tau_{min}+\tau_{max})/2$.
\item
If $(\tau_{max}-\tau_{min})/|\tau_{max}|>\varepsilon_{1}$, then go to 4).
Otherwise, $\tau\lto\tau_{max}$.
\item
If $\dstyle{\left|\frac{\partial F }{\partial\tau}(\lambda,\tau)\right|
<\varepsilon_{2}}$, then
record $\bigl(\lambda,\tau,F(\lambda,\tau)\bigr)$,
$\lambda\lto\lambda+\del{\lambda}$ and go to 2). 
Otherwise, end.
\end{enumerate}

In $\lambda<0$, replace 1), 2) and 6) with
$1^\prime$), $2^\prime$) and $6^\prime$), respectively:
\begin{enumerate}
\renewcommand{\labelenumi}{\arabic{enumi}$^\prime$)}
\item
$\tau\lto\delta_{0}$, $\lambda\lto-\varepsilon_{0}$.
\item
If $\lambda<\lambda_{min}^{-}$, end. Otherwise,
$\tau_{min}\lto -\lambda x_{1}$, $\tau_{max}\lto \tau_{U}^{-}(\lambda)$.
\addtocounter{enumi}{3}
\item
If $\dstyle{\left|\frac{\partial F }{\partial\tau}(\lambda,\tau)\right|
<\varepsilon_{2}}$, then
record $\bigl(\lambda,\tau,F(\lambda,\tau)\bigr)$,
$\lambda\lto\lambda-\del{\lambda}$ and go to 2$^\prime$). 
Otherwise, end.
\end{enumerate}
Here, $\delta_{0}$, $\varepsilon_{0}$, $\varepsilon_{1}$,
$\varepsilon_{2}$, $\lambda_{max}^{+}$, $\del{\lambda}$
and $\lambda_{min}^{-}$ are
preassigned constants for the procedure.
The way of determining them will be
explained in the next section.

Using $\bigl\{(\lambda,F(\lambda,\tau))\bigr\}$ in the records in 6)
and $6^{\prime}$),
we can plot the profile of $F(\lambda,\tau)$. In addition, if a PRM
exists and we set $\del{\lambda}$ at a sufficiently small positive value,
we can immediately get the extreme point
of $F(\lambda,\tau)$ with high accuracy.

\vskip 3mm

\noindent {\bf 4  COMPUTATIONAL EXPERIMENTS}

In this section we give searching examples and Monte Carlo studies. The
examples include the three types of cases: 1) the $\lambda$ coordinate of
a PRM is positive, 2) that is negative, 3) no PRM exists. The
Monte Carlo studies give a correlation among the value of the population
parameter $\lambda$, the rate at that $\lambda$ is positively
estimated and the existence rate of a PRM.
%

\vskip 3mm

\noindent {\bf 4.1  Searching Examples}

Using the algorithm in Section 3, we seek PRMs and profiles of
$F(\lambda,\tau)$, that is, $F(\lambda,\tau_{0}(\lambda))$ for six
data sets. They are indicated in Table I.
Data 1 has been introduced as an example to fail to find the PRM
in (Cohen, Whitten and Ding, 1985). Data 2 and 3 have been introduced
as difficult examples to seek the PRMs in (Cheng and Iles, 1990).
Data 4 and 5 are examples in which
$\hat{L}(\alpha)$ attains its maximum at a point
closely near to the end-point of the domain of definition.
We picked up Data 3 and 4
for the introduction in Section 1.
Data 6 is an artificial data set generated in Monte Carlo simulation
when the true values of $\lambda$, $\mu$ and $\sigma$ are
set at $0.4$, $0$ and $1$, respectively.

\input{tab1}

At first
the preassigned constants of the algorithm are selected in the
following way.
We set $\del{\lambda}=0.05$ so as to obtain a crude range in which
the $\lambda$ coordinate of
a PRM lies. 
Since zero is not included in the domain
of definition of the parameter $\lambda$, we set $\varepsilon_{0}
=0.05\ (>0)$ to exclude it from the search range.
Then, it is appropriate to set $\delta_{0}=\tau^{\ast}$
because of 5) in the theorem.
As shown in the last part of Section 2, the search range of $\lambda$
is sufficiently wide when $\lambda_{max}^{+}$ and
$\lambda_{min}^{-}$ are set at 6 and $-6$, respectively.
For each
$\lambda$, the $\tau_{0}(\lambda)$ can be calculated to an
accuracy of the order of the machine epsilon on a used computer. But
even if the calculation is performed with such
an accuracy, $|\frac{\partial F}{\partial \tau}(\lambda,\tau)|$ may
not take a value sufficiently close to $0$ for
a large $|\lambda|$ and
an approximate value of $\tau_{0}(\lambda)$. That is because
$\frac{\partial F}{\partial \tau}(\lambda,\tau)$ varies considerably
in $(-\lambda x_{n},\tau_{U}^{+}(\lambda))$ or
$(-\lambda x_{1},\tau_{U}^{-}(\lambda))$ when $|\lambda|$ is
large. Thus, we set that $\varepsilon_{1}=10^{-14}$ and
$\varepsilon_{2}=0.01$.

Under this setting, we have obtained the profiles of $F$.
Each one has a shape drawn with a solid curve on Fig. 3.
Dotted curves express $F(\lambda,\tau_{U}^{+}(\lambda))$
or $F(\lambda,\tau_{U}^{-}(\lambda))$. In the figure we can see
that solid and dotted curves almost overlap each other in the
intervals from $-2$ or $2$ up to
about $-4$ or $4$ in $\lambda$, respectively.

\clearpage
\begin{figure}[h]
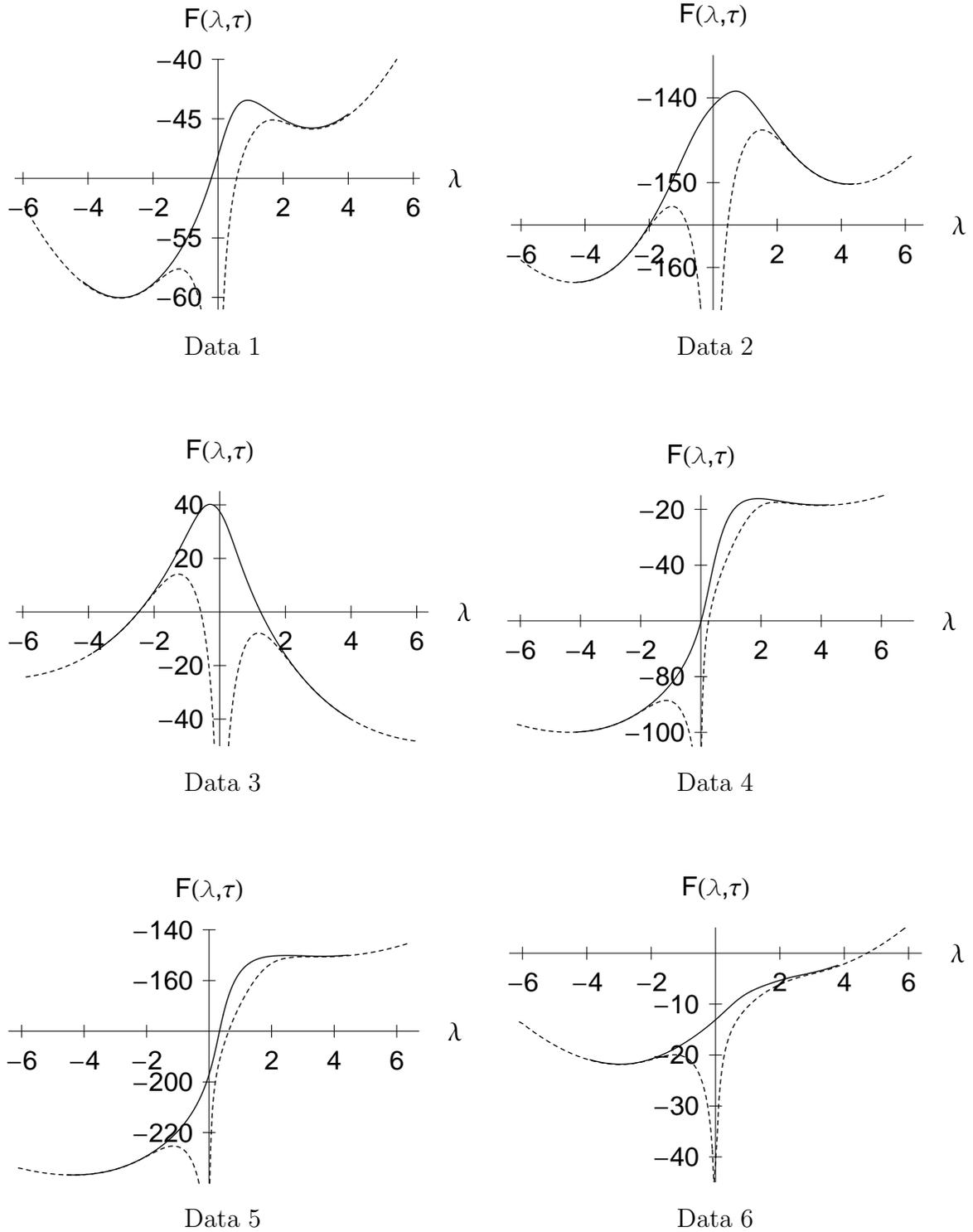

\unitlength 5mm
\begin{center}
\begin{picture}(32,40)(0,0)
\put(1,30){\epsfig{file=profile1.epsi,scale=1}}
\put(17,30){\epsfig{file=profile2.epsi,scale=1}}
\put(1,16){\epsfig{file=profile3.epsi,scale=1}}
\put(17,16){\epsfig{file=profile4.epsi,scale=1}}
\put(1,2){\epsfig{file=profile5.epsi,scale=1}}
\put(17,2){\epsfig{file=profile6.epsi,scale=1}}
\put(6.7,28.6){Data 1}
\put(22.5,28.6){Data 2}
\put(6.7,14.6){Data 3}
\put(22.5,14.6){Data 4}
\put(6.7,0.6){Data 5}
\put(22.5,0.6){Data 6}
\end{picture}
\caption{Profiles of $F(\lambda,\tau)$. Solid curves indicate
$F(\lambda,\tau_{0}(\lambda))$. Dotted curves indicate
$F(\lambda,\tau_{U}^{+}(\lambda))$ or $F(\lambda,\tau_{U}^{-}(\lambda))$.
These are included for comparison.}
\end{center}
\end{figure}
\clearpage

The profile for Data 1 indicates that the PRM exists and
its $\lambda$ coordinate is around $1$. In fact, we can
know that the $\lambda$ coordinate of the PRM is
in $(0.85,0.95)$ and $\tau_{0}(0.85)=-121.0106$
from the record

\vspace*{-8mm}
\begin{eqnarray*}
\{\bigl(\lambda,\tau,F(\lambda,\tau)\bigr)\}
&=&\{\ldots, (.85,-121.0106,-43.4512),
(.90,-129.4756,-43.4380),\\
& & \makebox[3em]{}
(.95,-137.8374,-43.4433),
\ldots
\}
\end{eqnarray*}
obtained by the algorithm.
Thus, if we perform only the part of the algorithm in $\lambda>0$ again
after resetting that
$\varepsilon_{0}=0.85$, $\delta_{0}=-121.0106$,
$\del{\lambda}=5\times 10^{-4}$ and
$\lambda_{max}^{+}=0.95$,
we can know the coordinate of the PRM more precisely,
which is $(\lambda,\tau)=(0.9095,-131.0716)$.
By means of a similar procedure, we obtain the results in Table II.
The profiles for  Data 2, 4 and 5 indicate a similar
situation to that for Data 1.
The profile for Data 3 differs from these
in the point that the $\lambda$ coordinate of the PRM is negative.
The last one is the profile for Data 6. It clearly shows that no PRM exists.
%

\input{tab2}
\vspace{5mm}


\vskip 3mm

\noindent {\bf 4.2  Monte Carlo Studies}

Varying the value of the parameter $\lambda$ in Monte Carlo simulation,
we investigate the existence rate of a PRM and
the rate at that $\lambda$ is positively estimated.
In addition, we seek the successful rate in
finding a PRM with Munro and Wixley's parameterization
for comparison's sake.

The simulation conditions are as follows: The sample number $n$ is
set at $10$, $15$ or $20$. The
parameter $\lambda$ is set at  $0.01$, $0.25$, $0.5$, $0.75,\ldots$,
$1.75$ or $2.0$, whereas the other parameters $\mu$ and $\sigma$ are
fixed at $0$ and $1$, respectively. For each combination of values
of $n$ and $\lambda$, $1000$ independent pseudo-random samples are 
considered.

We judge whether a PRM exists or not for each data set by tracing
$\{(\lambda, F(\lambda,\tau))\}$ in each record.
%

\input{tab3}

\input{tab4}

The simulation results for our parameterization are
shown in Table III and IV. From Table III we
can see that the existence rate of a PRM becomes lower as $\lambda$ becomes
larger. On the other hand, from Table IV we can see that the rate at
that $\lambda$ is positively estimated becomes higher as $\lambda$ becomes
larger.

\input{tab5}

The result for Munro and Wixley's parameterization is
shown in Table V. In the parameterization, even when
$\lambda$ is set at a value, it needs to solve
a two-dimensional non-linear equation not a scalar
equation. For it we used
Newton's method. Besides, we adopted
a similar algorithm to the one in Section 2.
That is, we set that
$\lambda\lto\pm 0.05$, $\mu\lto\bar{x}$ and
$\sigma\lto\sqrt{(1/n)\sum_{i=1}^{n}(x_{i}-\bar{x})^{2}}$
at the first stage (this corresponds to 1) or $1^{\prime}$)), and
chose, as the initial values of $\mu$ and $\sigma$ at the
present stage, the values of them obtained at the previous
stage (this corresponds to 3)).
Comparing Table III and V, we can see that
the failure rate in finding a PRM becomes higher
as $\lambda$ becomes larger from $0.5$.
In addition, these tables indicate that, although the existence rate is low
when $\lambda$ is large in the case of $n=10$, it is easier
to find a PRM than the other cases if it exists.

\vskip 3mm

\noindent {\bf 5   Discussion}

We discuss about $F(\lambda,\tau_{0}(\lambda))$ when $|\lambda|$ is
large. For each $\lambda$, $\tau_{0}(\lambda)$ lies in
the interval $(-\lambda x_{n}, \tau_{U}^{+}(\lambda))$
or $(-\lambda x_{1}, \tau_{U}^{-}(\lambda))$, and
as $|\lambda|\to\infty$ the widths of the intervals approach $0$
at a speed of $O(\e^{-\lambda^{2}+\ln |\lambda|})$.
Thus, $\tau_{0}(\lambda)$ approaches $-\lambda x_{n}$ or
$-\lambda x_{1}$ at a speed equal to or faster
than the order. Our question is how $F(\lambda,\tau_{0}(\lambda))$
behaves then.
If $|\frac{\partial F}{\partial \tau}(\lambda,\tau)|$ does
not increase rapidly for any
$\tau$ in $(\tau_{0}(\lambda),\tau_{U}^{+}(\lambda))$
or $(\tau_{0}(\lambda),\tau_{U}^{-}(\lambda))$
as $|\lambda|\to\infty$,
we can tell, with the mean value theorem,
that $F(\lambda,\tau_{0}(\lambda))$ behaves similarly to
$F(\lambda,\tau_{U}^{+}(\lambda))$
or $F(\lambda,\tau_{U}^{-}(\lambda))$ for large $|\lambda|$.
However, the $|\frac{\partial F}{\partial \tau}(\lambda,\tau)|$
may increase rapidly because
$|\frac{\partial F}{\partial \tau}(\lambda,\tau_{U}^{+}(\lambda))|$
and $|\frac{\partial F}{\partial \tau}(\lambda,\tau_{U}^{-}(\lambda))|$
increase at a speed of $O(\e^{\lambda^{2}-\ln |\lambda|})$.

We set that
$\tau^{+}=-\lambda x_{n}+g(\lambda)$ for $\lambda>0$ and
$\tau^{-}=-\lambda x_{1}+g(\lambda)$ for $\lambda<0$, where
$g(\lambda)>0$ and $g(\lambda)\to 0\ (\lambda\to\pm\infty)$. Then,
the following holds. As $|\lambda|\to\infty$
\begin{enumerate}
\renewcommand{\labelenumi}{\arabic{enumi}$)$}
\item
$F(\lambda, \tau^{+})\to -\infty$, $F(\lambda, \tau^{-})\to -\infty$ if
$g(\lambda) = o(\e^{-|\lambda|^{\eta}})$ for $\eta > 2$,
\item
$F(\lambda, \tau^{+})\to -\infty$, $F(\lambda, \tau^{-})\to -\infty$ if
$g(\lambda) = C_{1}\e^{-C_{2}\lambda^{2}+o(\lambda^{2})}$
for $C_{1}>0$, $C_{2}\geq\frac{2n}{n-1}$,
\item
$F(\lambda, \tau^{+})\to \infty$, $F(\lambda, \tau^{-})\to \infty$ if
$g(\lambda) = C_{1}\e^{-C_{2}\lambda^{2}+o(\lambda^{2})}$
for $C_{1}>0$, $0<C_{2}<\frac{2n}{n-1}$.
\end{enumerate}
This can be proved by standard mathematical calculations.

From 3), $\dstyle{\lim_{\lambda\to\pm\infty} F(\lambda,\tau_{0}
(\lambda)) = \infty}$ and the fact that $\tau_{0}(\lambda)$ lies
in $(-\lambda x_{n}, \tau_{U}^{+}(\lambda))$ or

{\noindent
$(-\lambda x_{1}, \tau_{U}^{-}(\lambda))$,
}
\[
\lambda x_{n}+\tau_{0}(\lambda)\sim
C_{1}\e^{-C_{2}^{\prime}\lambda^{2}+o(\lambda^{2})}
\quad {\rm or}\quad
\lambda x_{1}+\tau_{0}(\lambda)\sim
C_{1}\e^{-C_{2}^{\prime}\lambda^{2}+o(\lambda^{2})}
\]
for large $|\lambda|$. Here, $1\leq C_{2}^{\prime}<\frac{2n}{n-1}$.
Consequently, $F(\lambda,\tau_{0}(\lambda))$ is considered to behave
similarly to $F(\lambda,\tau_{U}^{+})$ or $F(\lambda,\tau_{U}^{-})$
for large $|\lambda|$. That is, we may regard
$F(\lambda,\tau_{U}^{+}(\lambda))$ or $F(\lambda,\tau_{U}^{-}(\lambda))$
as $F(\lambda,\tau_{0}(\lambda))$ for large $|\lambda|$.
See also the examples on Fig. 3.
\vskip 3mm

\noindent {\bf 6   SUMMARY}

We have proposed a reparameterization of the extended lognormal
distribution for the parameter estimation. The reparameterization
changes the three-parameter estimation problem to a two-parameter
estimation problem, and enables us to cope with extensive
data sets including
those which cause the embedded problem. On the two-parameter
estimation problem, we have made
an algorithm to seek the
profile of an object function. The algorithm is simple and makes it
possible to seek the profile stably. The profile clearly
shows whether a PRM exists or not, and if it exists, we can
obtain it certainly from the record for drawing the profile.
In fact, we have illustrated that the reparameterization and the algorithm
go well for the six data sets including ones introduced as difficult
examples to estimate the PRMs in other articles.

By means of Monte Carlo simulation, we have investigated the existence
rate of a PRM while varying the value of the parameter $\lambda$.
From the simulation result it has become clear
that the rate largely falls down when $\lambda$ goes through
from $1.25$ to $2.0$.
In such cases that the degree of skewness is high and no MLE
exists, there is a possibility that dealing with data as grouped data
can help us cope with the difficulty (Giesbrecht and Kempthorne, 1976).
This possibility was not pursued in the present paper because it
is beyond its purpose.
We have investigated the rate at that $\lambda$
is positively estimated on data governed by (1.2). The simulation result
has shown the rate is less than 1 in almost all the sample number
when $\lambda$ is less than 1.
This indicates the necessity for our or Munro and Wixley's
parameterization, which permits $\lambda$ to be negative,
since this generalization makes it possible to cope with the embedded problem.
\vskip 3mm

\noindent ACKNOWLEDGEMENTS
\vskip 3mm
The authors would like to thank the referees for their
helpful comments to improve this paper.
\vskip 3mm
\noindent {\bf\it References}
\vskip 3mm

\noindent Bilikan, J.E., Moore, A.H. and Petrick, G.L. (1979).
$K$ Sample ML ratio test for change of shape parameter,
{\it IEEE Trans. Reliab.}, {\bf 28}, 47--50.
\vskip 3mm

\noindent Chen, G. and Balakrishnan, N. (1995).
A general purpose approximate goodness-of-fit test,
{\it J. Quality Technol.}, {\bf 27} (2), 154--161.
\vskip 3mm

\noindent Cheng, R.C.H. and Iles, T.C. (1990).
Embedded models in three-parameter distributions and their estimation,
{\it J. Royal Statist. Soc.} B, {\bf 52} (1), 135--149.
\vskip 3mm

\noindent Cohen, A.C., Whitten, B.J. and Ding, Y. (1985).
Modified moment estimation for the three-parameter lognormal distribution,
{\it J. Quality Technol.}, {\bf 17} (2), 92--99.
\vskip 3mm

\noindent Eastham, J.F., LaRiccia, V.N. and Schuenemeyer, J.H. (1987).
Small sample properties of the maximum likelihood estimators for
an alternative parameterization of the three-parameter lognormal distribution,
{\it Comm. Statist.--Simulation Comput.}, {\bf 16} (3), 871--884.
\vskip 3mm

\noindent Giesbrecht, F. and Kempthorne, O. (1976).
Maximum likelihood estimation in the three-parameter lognormal distribution,
{\it J. Royal Statist. Soc.} B, {\bf 38}, 257--264.
\vskip 3mm

\noindent Hill, B.M. (1963).
The three-parameter lognormal distribution and
Bayesian analysis of a point-source epidemic,
{\it J. Amer. Statist. Assoc.}, {\bf 58}, 72--84.
\vskip 3mm

\noindent Hirose, H. (1997).
Maximum likelihood parameter estimation in the three-parameter log-normal
distribution using the continuation method,
{\it Comput. Statist. Data Anal}, {\bf 24}, 139--152.
\vskip 3mm

\noindent Johnson, N.L., Kotz, S. and Balakrishnan, N. (1994).
{\it Continuous Univariate Distributions}, New York, John Wiley \& Sons.
\vskip 3mm

\noindent Lambert, J.A. (1964).
Estimation of parameters in the three parameter lognormal distribution,
{\it Austral. J. Statist.}, {\bf 6}, 29--32.
\vskip 3mm

\noindent McCool, J.I. (1974).
Inferential techniques for Weibull populations,
{\it Aerospace Research Laboratories Report ARL TR 74-0180},
Wright-Patterson AFB, OH.
\vskip 3mm

\noindent Menon, M.V. (1963).
Estimation of the shape and scale parameters of
the Weibull distribution,
{\it Technometrics}, {\bf 5} (2), 175--182.
\vskip 3mm

\noindent Munro, A.H. and Wixley, R.A.J. (1970).
Estimation on order statistics
of small samples from a three-parameter lognormal distribution,
{\it J. Amer. Statist. Assoc.}, {\bf 65} (329), 212--225.
\vskip 3mm

\noindent Smith, R.L. and Naylor, J.C. (1987).
A comparison of maximum likelihood and Bayesian estimators for the
three-parameter Weibull distribution,
{\it Appl. Statist.}, {\bf 36}, 358--369.
\vskip 3mm

\noindent Steen, P.J. and Stickler, D.J. (1976).
A sewage pollution study of beaches from Cardiff to Ogmore,
UWIST, Dept. of Applied Biology Report, Cardiff.
\vskip 3mm

\noindent Wingo, D.R. (1975).
The use of interior penalty functions to overcome
lognormal distribution parameter estimation anomalies,
{\it J. Statist. Comput. Simulation}, {\bf 4}, 49--61.
\vskip 3mm

\noindent Wingo, D.R. (1976).
Moving truncations barrier-function methods for
estimation in three-parameter lognormal models,
{\it Comm. Statist.--Simulation Comput.}, {\bf B5} (1), 65--80.
\vskip 3mm

\noindent Wingo, D.R. (1984).
Fitting three-parameter lognormal models by
numerical global optimization--an improved algorithm,
{\it Comput. Statist. Data Anal.}, {\bf 2}, 13--25.
\vskip 3mm

\noindent APPENDIX
\vskip 3mm

The following is the first and second derivatives of the
object function $F$ with respect to $\tau$:
$$
\begin{array}[b]{l}
\dstyle{
\frac{\partial F}{\partial \tau}(\lambda,\tau)=
\frac{1}{n\lambda^2}
\left\{
\sum_{i=1}^{n}\ln(\lambda x_{i}+\tau)
\right\}
\left\{
\sum_{j=1}^{n}\frac{1}{\lambda x_{j}+\tau}
\right\}
}\\
%
\makebox[10em]{}
\raisebox{0mm}[7mm][0mm]{}
\dstyle{
-\frac{1}{\lambda^2}
\sum_{i=1}^{n}\frac{1}{\lambda x_{i}+\tau}\ln(\lambda x_{i}+\tau)
-\sum_{i=1}^{n}\frac{1}{\lambda x_{i}+\tau},
}
\end{array}
\eqno(A.1)
$$
$$
\begin{array}[b]{l}
\dstyle{
\frac{\partial^2 F}{\partial \tau^2}(\lambda,\tau)=
\frac{1}{n\lambda^2}
\left[
\left(\sum_{i=1}^{n}\frac{1}{\lambda x_{i}+\tau}\right)^2
-\left\{\sum_{j=1}^{n}\ln(\lambda x_{j}+\tau)\right\}
\left\{\sum_{j=1}^{n}\frac{1}{(\lambda x_{j}+\tau)^2}\right\}
\right.
}
\\
\makebox[8em]{}
\raisebox{0mm}[7mm][0mm]{}
\dstyle{
\left.
+\sum_{i=1}^{n}\frac{n}{(\lambda x_{i}+\tau)^2}\ln(\lambda x_{i}+\tau)
-\sum_{i=1}^{n}\frac{n}{(\lambda x_{i}+\tau)^2}
\right]
+\sum_{i=1}^{n}\frac{1}{(\lambda x_{i}+\tau)^2}.
}
\end{array}
\eqno(A.2)
$$
\end{document}

%% file: tab1.tex
\begin{table}[t]
\begin{center}
TABLE I \hspace{5mm} Data sets.
\vspace*{2mm}
\renewcommand{\arraystretch}{0.6}
\begin{tabular}{ccccccccccc}
\hline
\multicolumn{11}{l}{Data 1: fatigue life in hours of 10
bearings (Cohen {\it et al.}, 1985; McCool, 1974)}\\
152.7 & 172.0 & 172.5 & 173.3 & 193.0 & 204.7 & 216.5 & 234.9 & 262.6 & 422.6\\
\multicolumn{11}{l}{}\\
\multicolumn{11}{l}{Data 2: times to failure of vehicles (Cheng and Iles,
1990; Bilikan, Moore and Petrick, 1979)}\\
184  &   250 &  439 &  444 &  450 &  478 &  487 &  524 &  688 &  850 \\
1048 &  1280 & 1364 & 1488 & 1513 & 1860 & 1947 & 1991 & 2200 & 2446 \\
\multicolumn{11}{l}{}\\
\multicolumn{11}{l}{Data 3: strengths of 15 cm
fibres (Cheng and Iles, 1990; Smith and Naylor, 1987)}\\
0.37 & 0.40 & 0.70 & 0.75 & 0.80 & 0.81 & 0.83 & 0.86 & 0.92 & 0.92 \\
0.94 & 0.95 & 0.98 & 1.03 & 1.06 & 1.06 & 1.08 & 1.09 & 1.10 & 1.10 \\
1.13 & 1.14 & 1.15 & 1.17 & 1.20 & 1.20 & 1.21 & 1.22 & 1.25 & 1.28 \\
1.28 & 1.29 & 1.29 & 1.30 & 1.35 & 1.35 & 1.37 & 1.37 & 1.38 & 1.40 \\
1.40 & 1.42 & 1.43 & 1.51 & 1.53 & 1.61 \\
\multicolumn{11}{l}{}\\
\multicolumn{11}{l}{Data 4: Menon's data example (Menon, 1963)}\\
\multicolumn{1}{c}{$\e^{-6.824}$} & \multicolumn{1}{c}{$\e^{-3.506}$} &
\multicolumn{1}{c}{$\e^{-2.64}$}  & \multicolumn{1}{c}{$\e^{-1.686}$} &
\multicolumn{1}{c}{$\e^{-1.064}$} & \multicolumn{1}{c}{$\e^{-0.832}$} &
\multicolumn{1}{c}{$\e^{-0.758}$} & \multicolumn{1}{c}{$\e^{-0.754}$} &
\multicolumn{1}{c}{$\e^{-0.684}$} & \multicolumn{1}{c}{$\e^{ -0.438}$}
\\
\multicolumn{1}{c}{$\e^{-0.41}$}  & \multicolumn{1}{c}{$\e^{-0.216}$} &
\multicolumn{1}{c}{$\e^{-0.03}$}  & \multicolumn{1}{c}{$\e^{0.032}$}  &
\multicolumn{1}{c}{$\e^{0.438}$}  & \multicolumn{1}{c}{$\e^{0.716}$}  &
\multicolumn{1}{c}{$\e^{1.262}$}  & \multicolumn{1}{c}{$\e^{1.954}$}  &
\multicolumn{1}{c}{$\e^{2.208}$}  & \multicolumn{1}{c}{$\e^{4.054}$}
\\
\multicolumn{11}{l}{}\\
\multicolumn{11}{l}{Data 5: pollution data (Chen and Balakrishnan, 1995; Steen and Stickler, 1976)}\\
109 & 111 & 154 & 200 & 282 & 327 & 336 & 482 & 718 & 900 \\
918 & 1045 & 1082 & 1345 &  1415 & 1918 & 2120 & 5900 & 6091 & 53600\\
\multicolumn{11}{l}{}\\
\multicolumn{11}{l}{Data 6: an artificial data set
generated in Monte Carlo simulation}\\
\multicolumn{2}{r}{$-0.912527$} & \multicolumn{2}{r}{$-0.905886$} &
\multicolumn{2}{r}{$-0.836045$} & \multicolumn{2}{r}{$-0.382619$} &
\multicolumn{2}{r}{$-0.319501$}
\\
\multicolumn{2}{r}{$0.030242$}  & \multicolumn{2}{r}{$0.326860$}  &
\multicolumn{2}{r}{$2.325620$}  & \multicolumn{2}{r}{$4.333967$}  &
\multicolumn{2}{r}{$5.663170$}
\\
\hline
\end{tabular}
\end{center}
\end{table}

%% file: tab2.tex
\begin{table}[h]
\begin{center}
TABLE II \hspace{5mm} Estimates of $\lambda$ and $\tau$.

\vspace*{-5mm}
\unitlength 1cm
\begin{picture}(16.5,2)
\put(2,0){
\begin{tabular}{c|ccccc}
\hline
 & \multicolumn{5}{|c}{Data} \\
 & $1$ & $2$ & $3$ & $4$ & $5$ \\
\hline
$\lambda$ & $.9095$     & $.7030$   & $-0.2955$ & $1.9065$ & $2.5135$    \\
$\tau$    & $-131.0716$ & $28.3203$ & $.5984$   & $.0126$  & $-272.6434$ \\
\hline
\end{tabular}
}
\end{picture}
\end{center}
\end{table}

%% file: tab3.tex
\begin{table}[b]
\begin{center}
TABLE III \hspace{5mm} Existence rate of a PRM.

\vspace*{2mm}
\begin{tabular}{cc|ccccccccc}
\hline
& & \multicolumn{9}{|c}{$\lambda$} \\
& & $.01$ & $.25$ & $.50$ & $.75$
  & $1.0$  & $1.25$  & $1.5$ & $1.75$
  & $2.0$
\\
\hline
          &10 & $.98$ & $.97$ & $.96$ & $.91$
              & $.83$ & $.71$ & $.56$ & $.40$
              & $.28$
\\
$n$       &15 & $1.0$ & $1.0$ & $1.0$ & $1.0$
              & $.98$ & $.96$ & $.87$ & $.73$
              & $.54$
\\
          &20 & $1.0$ & $1.0$ & $1.0$ & $1.0$
              & $1.0$ & $1.0$ & $.98$ & $.93$
              & $.80$
\\
%
%
\hline
\end{tabular}
\end{center}
\end{table}

%% file: tab4.tex
\begin{table}[t]
\begin{center}
TABLE IV \hspace{5mm} Rate at that $\lambda$ is positively estimated.

\vspace*{2mm}
\begin{tabular}{cc|ccccccccc}
\hline
& & \multicolumn{9}{|c}{$\lambda$} \\
& & $.01$ & $.25$ & $.50$ & $.75$
  & $1.0$  & $1.25$  & $1.5$ & $1.75$
  & $2.0$
\\
\hline
          &10 & $.49$ & $.72$ & $.85$ & $.94$
              & $.97$ & $.99$ & $.98$ & $.99$
              & $.99$
\\
$n$       &15 & $.51$ & $.79$ & $.93$ & $.99$
              & $1.0$ & $1.0$ & $1.0$ & $1.0$
              & $1.0$
\\
          &20 & $.51$ & $.85$ & $.98$ & $1.0$
              & $1.0$ & $1.0$ & $1.0$ & $1.0$
              & $1.0$
\\
%
%
\hline
\end{tabular}
\end{center}
\end{table}

%% file: tab5.tex
\begin{table}[t]
\begin{center}
TABLE V \hspace{5mm} Successful rate in finding a PRM with Munro and
Wixley's parameterization.

\vspace*{2mm}
\begin{tabular}{cc|ccccccccc}
\hline
& & \multicolumn{9}{|c}{$\lambda$} \\
& & $.01$ & $.25$ & $.50$ & $.75$
  & $1.0$  & $1.25$  & $1.5$ & $1.75$
  & $2.0$
\\
\hline
          &10 & $.98$ & $.96$ & $.94$ & $.87$
              & $.74$ & $.58$ & $.42$ & $.27$
              & $.17$
\\
$n$       &15 & $1.0$ & $1.0$ & $.99$ & $.97$
              & $.85$ & $.67$ & $.42$ & $.24$
              & $.13$
\\
          &20 & $1.0$ & $1.0$ & $1.0$ & $.99$
              & $.91$ & $.69$ & $.39$ & $.18$
              & $.08$
\\
\hline
\end{tabular}
\end{center}
\end{table}